\documentclass[10pt]{article}
\usepackage{amsfonts}
\usepackage[dvips]{graphicx}
\usepackage{amssymb,color}
\usepackage[latin1]{inputenc}
\usepackage{epic}
\usepackage[latin1]{inputenc}
\usepackage{enumerate}
\usepackage{color}
\usepackage{lettrine}
\usepackage{calc,pifont}
\usepackage[noindentafter,calcwidth]{titlesec}
\usepackage{amsmath}
\usepackage{amssymb}
\usepackage{amsthm}
\usepackage{subfigure}
\usepackage{lineno}

\newtheorem{lemma}{Lemma}[section]

\newtheorem{theorem}{Theorem}[section]

\newtheorem{corollary}{Corollary}[section]
\newtheorem{proposition}{Proposition}[section]

\theoremstyle{remark}
\newtheorem{rem}{\bf Remark}[section]
\theoremstyle{remark}

\title{Domination and 2-degree-packing numbers in graphs}

\author{Adri{\' a}n V\'azquez-\'Avila  \footnotemark[1]}
\date{}

\begin{document}
\maketitle

\def\thefootnote{\fnsymbol{footnote}}
\footnotetext[1]{Subdirecci{\' o}n de Ingenier{\' i}a y Posgrado, Universidad Aeronáutica en Quer\'etaro, Parque Aeroespacial de Quer\'etaro, 76278, Quer\'etaro, M\'exico, {\tt adrian.vazquez@unaq.edu.mx}.}
\vspace{-.5cm}
\begin{abstract}
A dominating set of a graph $G$ is a set $D\subseteq V(G)$ such that \-every vertex of $G$ is either in $D$ or is adjacent to a vertex in $D$. The domination number of $G$, $\gamma(G)$, is the minimum order of a domi\-nating set. A subset $R$ of edges of a graph $G$ is a 2-degree-packing, if any three edges from $R$ do not have the same incident vertex. The 2-degree-packing number of $G$, $\nu_2(G)$, is the maximum order of a 2-degree-packing of $G$.

In this paper, we prove that any simple graph $G$ satisfies $\gamma(G)\leq\nu_2(G)-1$. Furthermore, we give a characterization of simple connected graphs $G$ satisfying $\gamma(G)=\nu_2(G)-1$. 
\end{abstract}

Keywords. Domination, covering, 2-degree-packing.

Math. Subj. Class.: 05C69, 05C70.
\section{Introduction}\label{sec1}

In this paper, we consider finite undirected simple graphs. For any undefined terms see \cite{Bondy}. 

Let $G$ be a graph with set of vertices $V(G)$ and set of edges $E(G)$. For a subset $A\subseteq V(G)$, $G[A]$ denotes the subgraph of $G$ which is \emph{induced by the vertex set} $A$. At the same way, for a subset $R\subseteq E(G)$, $G[R]$ denotes the subgraph of $G$ which is \emph{induced by the edge set} $R$. The \emph{open neighborhood} of a vertex $u\in V(G)$, $N(u)$, is the set of vertices of $V(G)$ adjacent to $u$ in $G$, and the \emph{closed neighborhood} of a vertex $u\in V(G)$, $N[u]$, is defined as $N[u]=N(u)\cup\{u\}$. The \emph{degree} of a vertex $u\in V(G)$, $deg(u)$, is defined as $deg(u)=|N(u)|$. The maximum degree of the graph $G$ and the minimum degree of the graph $G$ are denoted by $\delta(G)$ and $\Delta(G)$, respectively. Let $H$ be a subgraph of $G$. The \emph{restricted open neighborhood} for a vertex $u\in V(H)$, $N_H(u)$, is defined as $N_H(u)=\{v\in V(H):uv\in E(H)\}$, the \emph{restricted closed neighborhood} for a vertex $u\in V(H)$, $N_H[u]$, is defined as $N_H[u]=N_H(u)\cup\{u\}$, and the \emph{restricted degree} of a vertex $u\in V(H)$, $deg_H(u)$, is defined as $deg_H(u)=|N_H(u)|$.

An \emph{independent set} of a graph $G$ is a subset $I\subseteq V(G)$ such that any two vertices of $I$ are not adjacent. The \emph{independence number} of $G$, $\alpha(G)$, is the maximum order of an independent set. A \emph{vertex cover} of a graph $G$ is a subset $T\subseteq V(G)$ such that all edges of $G$ has at least one end in $T$. The \emph{covering number} of $G$, $\beta(G)$, is the minimum order of a vertex cover of $G$. A \emph{dominating set} of a graph $G$ is a set $D\subseteq V(G)$ such that each vertex $u\in V(G)\setminus D$ satisfies $N(u)\cap D\neq\emptyset$. The \emph{domination number} of $G$, $\gamma(G)$, is the minimum order of a dominating set.  
 
It is well-known, if $G$ is a graph with no isolated vertices, then
\begin{equation}\label{eq:gamma}
\gamma(G)\leq\beta(G).
\end{equation}

In \cite{RV98}, was given a characterization of simple graphs which attains the inequality (\ref{eq:gamma})  (see also \cite{VOLKMANN1994211,Yun}).

A \emph{k-degree-packing set} of a graph $G$ ($k\leq\Delta(G)$), is a subset $R\subseteq E(G)$ such that $\Delta(G[R])\leq k$. The \emph{k-degree-packing number} of $G$, $\nu_k(G)$, is the maximum order of a $k$-degree-packing set of $G$, see \cite{CACA}. We are interested when $k=2$, since $k=1$ is the \emph{matching number} of a graph. 

The 2-degree-packing number is studied in \cite{CGCA,AGAL,Avila,Avila2,Avila3,CACA} in a more general context, but with a different name, as 2-packing number, see \cite{CACA}. A 2-packing in graphs has different meaning: A set $X\subseteq V(G)$ is called a 2-packing if $d_G(u,v)>2$ for any different vertices $u$ and $v$ of $X$, that is, the2-packing is a subset $X\subseteq V(G)$ in which all the vertices are in distance at least 3 from each other, see for example \cite{TOPP1991229}.
 
In \cite{CACA}, was proved the following chain of inequalities for a simple connected graph $G$:
\begin{equation}\label{eq:cubierta}
\lceil \nu_{2}(G)/2\rceil\leq\beta(G)\leq\nu_2(G)-1.
\end{equation}

Moreover, in \cite{CACA}, was given a characterization of simple connected graphs $G$ which attain the upper and lower inequality of (\ref{eq:cubierta}). 

Hence, by the chain of inequalities (\ref{eq:gamma}) and (\ref{eq:cubierta}), we have 
\begin{equation}\label{eq:des}
\gamma(G)\leq\nu_2(G)-1.
\end{equation}

In this paper, we give a characterization of simple connected graphs $G$ satisfying $\gamma(G)=\nu_2(G)-1$. There is not a lower bound for the domination number of a graph in terms of the 2-degree-packing number, since  $\gamma(K_n)=1$ and $\nu_2(K_n)=n$, for all $n\geq3$.  
\section{Main result}
To begin with, we introduce some terminology in order to simplify the description of the simple connected graphs $G$ which satisfy $\gamma(G)=\nu_2(G)-1$.

\begin{figure}[t]
	\begin{center}
		\includegraphics[height =2.7cm]{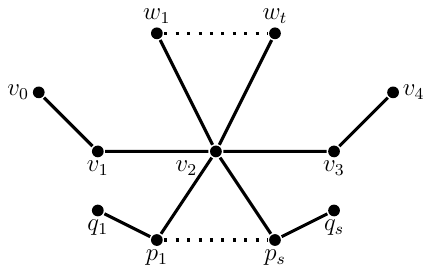}
	\end{center}
	\caption{$T^{r}_{s,t}$,}
	\label{fig:arbol_gamma}
\end{figure}

Let $P_4$ be a path of length 4, say $P_4=v_0v_1v_2v_3v_4$. We define the tree $T_{s,t}^{r}=(V,E)$, with $s+4=r$, as follow:
\begin{eqnarray*}
V&=&V(P_4)\cup\{p_1,\ldots,p_{s}\}\cup\{q_1,\ldots,q_{s}\}\cup\{w_1,\ldots,w_t\},\\
E&=&E(P_4)\cup\{p_iq_i:i=1,\ldots,s\}\cup\{v_2w_i:i=1,\ldots,t\}\cup\{v_2p_i:i=1,\ldots,s\}.
\end{eqnarray*}
and depicted in Figure \ref{fig:arbol_gamma}.

\begin{proposition}\label{prop:gamma}
	$\gamma(T_{s,t}^{r})=r-1$.
\end{proposition}
\begin{proof}
It is not difficult to see that $D=(V(P_4)\setminus\{v_0,v_4\})\cup\{p_1,\ldots,p_s\}$ is a dominating set of $T_{s,t}^{r}$ of cardinality $s+3=r-1$, which implies that $r-1\geq\gamma(T_{s,t}^{r})$. On the other hand, let $D\subseteq V(T_{s,t}^{r})$ be a minimum dominating set of $T_{s,t}^{r}$ with the minimum number of vertices of degree one. Without loss of generality, we assume that $deg(u)\geq2$, for every $u\in D$, otherwise $(D\setminus\{u\})\cup N(u)$ is a dominating set with less vertices of degree one, which is a contradiction. Hence, $\gamma(T_{s,t}^{r})=|D|\geq s+3=r-1$. Therefore, $\gamma(T_{s,t}^{r})=r-1$, and the statement holds.	
\end{proof}

\begin{proposition}\label{prop:2-packing}
	$\nu_2(T_{s,t}^r)=r$.
\end{proposition}

\begin{proof}
It is not difficult to see that $R=E(P_4)\cup\{p_iq_i:i=1,\ldots,s\}$
is a 2-degree-packing of $T_{s,t}^{r}$ of cardinality $s+4=r$. Hence,  $\nu_2(T_{s,t}^{r})\geq r$. 

On the other hand, notice that $deg(v_2)=r+t-2$ while $|E(T_{s,t}^r)|=2r+t-4$. Hence, the number of adges non-incident to $v_2$ is $r-2$, which implies that any set of $r+1$ edges must contain at least three edges incident to $v_2$, implying that any set of $r+1$ edges may not be a 2-degree-packing of $T_{s,t}^r$. Hence, $\nu_2(T_{s,t}^r)\leq r$.
\end{proof}

Hence, we have the following:

\begin{corollary}\label{coro:vuelta}
$\gamma(T_{s,t}^r)=\nu_2(T_{s,t}^r)-1$.		
\end{corollary}

The main result of this paper is stated as follow
\begin{theorem}
If $G$ is a simple connected graph with $\gamma(G)=\nu_2(G)-1$, then $G\simeq T^{\nu_2}_{s,t}$, where $\nu_2=\nu_2(G)$.	
\end{theorem}

\begin{figure}[t]
	\begin{center}
		\subfigure[]{\includegraphics[height=2.5cm]{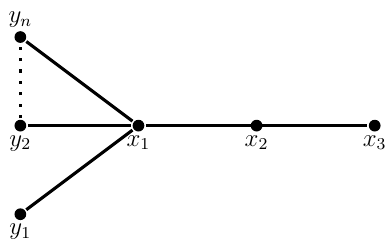}}
		\hspace{1 cm}
		\subfigure[]{\includegraphics[height=1.8cm]{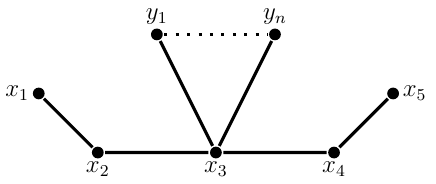}}
	\end{center}
	\caption{In $(a)$ the only family of graphs $G$ which satisfies $\nu_2(G)=3$ and $\gamma(G)=2$ is shown. On the other hand, in $(b)$ the only family of graphs which satisfies $\nu_2(G)=4$ and $\gamma(G)=3$ is shown.}
	\label{fig:nu_2_3}
\end{figure}

From here in follow, we are going to consider simple connected graphs with $|E(G)|>\nu_2(G)$, due to the fact $|E(G)|=\nu_2(G)$, if and only if, $\Delta(G)\leq2$.  On the other hand, since $\gamma(G)\leq\beta(G)$, for a simple connected graph $G$, then the bipartite graph $K_{1,n}$, for $n\geq3$, is the unique graph with $\gamma(G)=1$ and $\nu_2(K_{1,n})=2$, since $\nu_2(G)=2$, if and only if, $\tau(G)=1$, see \cite{AGAL}. It is not difficult to prove, the families of graphs of Figure \ref{fig:nu_2_3} are the unique families of graphs $G$ which satisfy $\nu_2(G)=3$ and $\beta(G)=2$, incise $(a)$, and $\nu_2(G)=4$ and $\beta(G)=3$, incise $(b)$, see \cite{AGAL,CACA}. Hence, we assume $\nu_2(G)\geq5$. 

\begin{lemma}\label{lemma:conected}
Let $G$ be a simple connected graph with $|E(G)|>\nu_2(G)$, and let $R$ be a maximum 2-degree-packing of $G$. If $G[R]$ is a connected graph, then $\gamma(G)\leq\nu_2(G)-2$. 
\end{lemma}
\begin{proof}
Let $R$ be a maximum 2-degree-packing of $G$ such that $G[R]$ is a connected graph, and let $I=V(G)\setminus V(G[R])$. Hence, if $I\neq\emptyset$, then $I$ is an independent set. 
	\begin{itemize}
		\item [Case (i):] $I=\emptyset$. Hence, $V(G)=V(G[R])=\{u_0,u_1,\ldots,u_{\nu_2(G)}\}$ and $|V(G)|=\nu_2(G)+1\geq6$, then $D=V(G)\setminus\{u_0,u_2,u_{\nu_2(G)}\}$ is a dominating set of size $\nu_2(G)-2$.

		\item [Case(ii):] $I\neq\emptyset$. It is well-known that either $G[R]$ is a path or a cycle. Let's suppose $G[R]$ is a path, say $G[R]=u_0u_1\cdots u_{\nu2}$. If $u\in I$, then $u_0,u_{\nu_2}\not\in N(u)$, otherwise either $R\cup\{u_0u\}$ or $R\cup\{uu_{\nu_2}\}$ is a 2-degree-packing of cardinality $\nu_2(G)+1$, which is a contradiction. On the other hand, if there is $u_i\in D^{**}=V(G[R])\setminus\{u_0,u_1,u_{\nu_2-1},u_{\nu_2}\}$ (since $\nu_2(G)\geq5$) such that $u_i\not\in N(u)$, for all $u\in I$, then $\hat{D}=V(G[R])\setminus\{u_0,u_i,u_{\nu_2}\}$ is a dominating set of cardinality $\nu_2(G)-2$, and the statement holds. Let's suppose that for every $u_i\in D^{**}$ there is $u\in I$ such that $u_i\in N(u)$. If $u_i\in N(u)$, for some $u\in I$, then $u_{i+1}\not\in N(u)$, otherwise the following set $\hat{R}=(R\setminus\{u_iu_{i+1}\})\cup\{uu_i,uu_{i+1}\}$ is a 2-degree-packing of $G$ of cardinality $\nu_2(G)+1$, which is a contradiction. Therefore, if $u_i,u_{i+1}\in D^{**}$, then there are $u,u^\prime\in I$ such that $u_i\in N(u)$ and $u_{i+1}\in N(u^\prime)$. Hence the following set $\hat{R}=(R\setminus\{u_iu_{i+1}\})\cup\{uu_i,u^\prime u_{i+1}\}$ is a 2-degree-packing of $G$ of cardinality $\nu_2(G)+1$, which is a contradiction. 
		
		To end the proof, let's suppose $G[R]$ is a cycle, say $G[R]=u_0u_1\cdots u_{\nu_2-1}u_0$. If there exists two different vertices $u_i,u_j\in V(G[R])$ such that $N(u_i)\cap I=\emptyset$ and $N(u_j)\cap I=\emptyset$, then the set $V(G[R])\setminus\{u_i,u_j\}$ is a dominating set of size $\nu_2(G)-2$. Alternatively, if at most one vertex in the cycle $G[R]$ is not adjacent to an element of $I$, then there exists adjacent vertices $u_i,u_{i+1}\in V(G[R])$ and $u,u^\prime\in I$ such that $u\in N(u_i)$ and $u^\prime\in N(u_j)$ (it is possible that $u=u^\prime$). Then the set   $$\hat{R}=(R\setminus\{u_iu_{i+1}\})\cup\{uu_i,u^\prime u_{i+1}\}$$ is a 2-degree-packing of size $\nu_2(G)+1$, contradiction. 
	\end{itemize} 	
\end{proof}

\begin{lemma}\label{lemma:forest}
Let $G$ be a simple connected graph with $|E(G)|>\nu_2(G)$ and $R$ be a maximum 2-degree-packing of $G$. If $\gamma(G)=\nu_2(G)-1$, then $G[R]$ is a forest. 	
\end{lemma}
\begin{proof}
Let $R$ be a maximum 2-degree-packing of $G$. The induced graph $G[R]$ is not a connected graph (by Lemma \ref{lemma:conected}). Hence, let $R_1, R_2,\ldots R_k$, with $k\leq\nu_2(G)-1$, be the components of $G[R]$, with $k$ as small as possible. Let's suppose $R_1,\dots,R_s$, are the components of $G[R]$ with only one edge (that is $R_i\simeq K_2$, for $i=1,\ldots,s$), that is $R_i=p_iq_i$, for $i=1,\ldots,s$, and $R_{s+1},\ldots,R_k$ are the components with at least two edges. Trivially $$D=\{u\in V(G[R]):deg_{R}(u)=2\}\cup\{p_i\in V(R_i):i=1,\ldots,s\}$$ is a dominating set of cardinality at most $\nu_2(G)$. 
	
Let $I=V(G)\setminus V(G[R])$, then
	
	\begin{itemize}
		\item [Case (i):] $I=\emptyset$. Let's suppose $R_{s+1}$ is a cycle and there is an edge $uv\in E(G)$ such that $u\in V(R_{s+1})$ and $v\in V(R_j)$, for some $j\in\{1,\ldots,k\}\setminus\{s+1\}$. Then $\hat{D}=D\setminus N_{R_{s+1}}(u)$ is a dominating set of $G$ of cardinality $\nu_2(G)-2$, a contradiction. Therefore, there are not cycles as components of $G[R]$.
		\item [Case (ii):] $I\neq\emptyset$. Let's suppose $R_{s+1}$ is a cycle and there is $u\in I$ such that $uv_{s+1}\in E(G)$, where $v_{s+1}\in V(R_{s+1})$.
		\begin{rem}\label{remark:1}
			If $v,w\in V(R_{s+1})$ are such $vw\in E(R_{s+1})$, then $v,w\not\in N(u)$, otherwise the following set $\hat{R}=(R\setminus\{vw\})\cup\{uv,uw\}$ is a 2-degree-packing of $G$ of cardinality $\nu_2(G)+1$, a contradiction. Simi\-larly,	if $v,w\in V(R_{s+1})$ are such $vw\in E(R_{s+1})$, then there are not $u,u^\prime\in I$ such that $v\in N(u)$ and $w\in N(u^\prime)$, otherwise the following set $\hat{R}=(R\setminus\{vw\})\cup\{uv,u^\prime w\}$ is a 2-degree-packing of $G$ of cardinality $\nu_2(G)+1$, which is a contradiction.
		\end{rem}
		By Remark \ref{remark:1}, we have $\hat{D}=D\setminus N_{R_{s+1}}(v_{s+1})$ is a dominating set of $G$ of cardinality $\nu_2(G)-2$, which is a contradiction. Therefore there are not cycles as components of $G[R]$.
	\end{itemize}
\end{proof}

\begin{theorem}\label{thm:ida}
Let $G$ be a simple connected graph with $|E(G)|>\nu_2(G)$. If $\gamma(G)=\nu_2(G)-1$, then $G\simeq T^{\nu_2}_{s,t}$, where $\nu_2=\nu_2(G)$. 
\end{theorem}
\begin{proof}
Let $R$ be a maximum 2-degree-packing of $G$ and $R_1,\ldots,R_k$ be the components of $G[R]$, with $2\leq k\leq\nu_2(G)-1$ (by Lemma \ref{lemma:conected}) and $k$ as small as possible. By Lemma \ref{lemma:forest} each component of $G[R]$ is a path. Let's suppose that $R_1,\dots,R_s$ are the components of $G$ with only one edge, that is, $R_i=p_iq_i$, for $i=1,\ldots,s$, and $R_{s+1},\ldots,R_k$ are the components with at least two edges. Trivially $$D=\{u\in V(G[R]):deg_{R}(u)=2\}\cup\{p_i\in V(R_i):i=1,\ldots,s\},$$ is a dominating set of cardinality at most $\nu_2(G)-1$. This implies either there is at most one component of length greater or equal than 2 and the rest of the components have only one edge, or all components of $G[R]$ have only one edge. Since $k$ is as small as possible, then $R_1,\ldots,R_{k-1}$ are the components with only one edge and $R_k$ is a path of length greater or equal than 2. 
	\begin{rem}\label{remark:2}
		If $u\in I=V(G)\setminus V(G[R])$, then neither $up_i\in E(G)$ nor $uq_i\in E(G)$, for all $i=1,\ldots,k-1$, otherwise either  $\hat{R}=R\cup\{up_i\}$ or $\hat{R}=R\cup\{uq_i\}$ is a 2-degree-packing of $G$ of size $\nu_2(G)+1$, a contradiction.	
	\end{rem}

Let's suppose that $R_k=v_0\cdots v_l$, with $l=2,3$. Since $\Delta>2$, then there is an edge $v_1w\in E(G)$ (possible $v_2w\in E(G)$ if $l=3$), with $w\in V(R_t)$, for some $t\in\{1,\ldots,k-1\}$. Hence, either $\hat{R}=(R\setminus\{v_0v_1\})\cup\{v_1w\}$ or $\hat{R}=(R\setminus\{v_2v_3\})\cup\{v_2,w\}$ is a maximum 2-degree-packing of $G$ with less components than $R$, which is a contradiction. Therefore, $R_k=v_0v_1\cdots v_l$, where $l\geq4$.
	
Let $V_k^{*}=V(R_k)\setminus\{v_0,u,v_l\}$, where $u$ is adjacent to either $u_0$ or $u_l$. If $v_ip_t\in E(G)$, with $v_i\in V^*_k$ and $p_t\in V(R_t)$, for some $t\in\{1,\ldots,k-1\}$, then $deg(q_t)=1$, otherwise either $v_iq_t\in E(G)$ or $v_j q_t\in E(G)$, for some $v_j\in V_k^*\setminus\{u_i\}$ and with $i\leq j$. In both cases either $\hat{R}=(R\setminus N_{R_k}(u_i))\cup\{u_ip_t,u_i q_t\}$ or  $\hat{R}=(R\setminus\{v_{i-1}v_i,v_jv_{j+1}\})\cup\{v_ip_t,v_jq_t\}$, is a maximum 2-degree-packing of $G$ with a cycle as a component, which is a contradiction (by Lemma \ref{lemma:forest}). Therefore, for all $i\in\{1,\ldots,k-1\}$ there is $v\in V_k^*$ such that $p_iv\in E(G)$ and $deg(q_i)=1$. 
	\begin{rem}\label{remark:3}
		If $u\in I$, then $uv\in E(G)$, for some $v\in V^*_k$, otherwise either $\hat{R}=R\cup\{v_0u\}$ is a 2-degree-packing of $G$ of cardinality $\nu_2(G)+1$, or $\hat{R}=(R\setminus\{v_1v_2\})\cup\{v_1u\}$ is a 2-degree-packing of $G$ with two components as paths of length at least 2, a contradiction. On the other hand, if either $uu_1\in E(G)$ or $uu_{l-1}\in E(G)$, then either $\hat{R}=(R\setminus\{u_0u_1\})\cup\{u_1u\}$ or $\hat{R}=(R\setminus\{u_{l-1}u_l\})\cup\{u_{l-1}u\}$ is a 2-degree-packing with less components, a contradiction. Also, if $uv\in E(G)$  with $v\in\{u_0,u_l\}$, then $\hat{R}=R\cup\{uv\}$ is a 2-degree packing of $G$ of size $\nu_2(G)+1$.
	\end{rem}

Now, we will prove that $l=4$, that is, $R_k=v_0\cdots v_4$. Suppose that $R_k=v_0\cdots v_l$, with $l\geq5$. If $v_i,p_j\in E(G)$, for some $v_i\in D\setminus\{u_1,u_{l-1}\}$ and $j\in\{1,\ldots,s\}$, then either $\hat{R}=(R\setminus\{v_iv_{i+1}\})\cup\{v_ip_j\}$ or $\hat{R}=(R\setminus\{v_{i-1}v_i\})\cup\{v_ip_j\}$ is a maximum 2-degree-packing of $G$ with two components as paths of length greater than 2, a contradiction. Therefore $l=4$, $v_2p_i\in E(G)$ and $v_2u\in E(G)$, for all $i=1,\ldots,s$ and $u\in I$. 
	
Finally, to show that $G\simeq T_{s,t}^{\nu_2}$, it remains to
be verified that $v_1v_3,v_0v_2,v_2v_4\not\in E(G)$ and $I\neq\emptyset$. If $v_1v_3\in E(G)$, then $\hat{R}=(R\setminus\{v_1v_2,v_3v_4\})\cup\{v_1v_3,v_2p_1,\}$ is a maximum 2-degree-packing with a path of length 5 as a component, which is a contradiction. On the other hand, if either $v_2v_4\in E(G)$ or $v_0v_2\in E(G)$, then either $\hat{R}=(R\setminus\{v_1v_2\})\cup\{v_2v_4\}$ or $\hat{R}=(R\setminus\{v_2v_3\})\cup\{v_0v_2\}$ is a 2-degree-packing which contains a cycle, a contradiction. Thus, the graph $G$ does not contain a cycle as a subgraph of $G$. To end, if $I=\emptyset$, then $D\setminus\{v_2\}$ is a dominating set of $G$ of cardinality $\nu_2(G)-2$, which is a contradiction. Therefore $G\simeq T_{s,t}^{\nu_2}$. 	
\end{proof}

Hence, the main results of this paper states:

\begin{theorem}
	Let $G$ be a connected graph with $|E(G)|>\nu_2(G)$. Then $\gamma(G)=\nu_2(G)-1$ if and only if $G\simeq T_{s,t}^{\nu_2}$.
\end{theorem}

\

{\bf Acknowledgment}

\


Research was partially supported by SNI, México and CONACyT, México.

\end{document}